\documentclass[11pt]{article}
\usepackage[T1]{fontenc}
\usepackage{lmodern,amsmath,amssymb,amsthm,mathtools,booktabs,array}
\usepackage[margin=1in]{geometry}
\usepackage{microtype}
\usepackage[hidelinks,hyperfootnotes=false]{hyperref}
\hypersetup{pdftitle={Long-root SL2 dual pairs in split exceptional groups},
 pdfauthor={Gordan Savin},
 pdfsubject={Compact types, cone geometry, and theta correspondence}}
\newcommand{\C}{\mathbb C}
\newcommand{\R}{\mathbb R}
\newcommand{\Z}{\mathbb Z}

\newcommand{\PP}{\mathbb P}
\newcommand{\OO}{\mathcal O}
\newcommand{\Gr}{\operatorname{Gr}}
\newcommand{\LG}{\operatorname{LG}}
\newcommand{\OGr}{\operatorname{OGr}}
\newcommand{\Sym}{\operatorname{Sym}}
\newcommand{\Hom}{\operatorname{Hom}}
\newcommand{\gr}{\operatorname{gr}}
\newcommand{\ch}{\operatorname{ch}}
\newcommand{\wt}{\operatorname{wt}}

\newcommand{\Irr}{\operatorname{Irr}}
\newcommand{\SL}{\operatorname{SL}}
\newcommand{\SO}{\operatorname{SO}}
\newcommand{\SU}{\operatorname{SU}}
\newcommand{\Sp}{\operatorname{Sp}}
\newcommand{\Spin}{\operatorname{Spin}}
\newcommand{\fg}{\mathfrak g}
\newcommand{\fk}{\mathfrak k}
\newcommand{\fp}{\mathfrak p}
\newcommand{\fh}{\mathfrak h}
\newcommand{\bS}{\mathbf S}
\newtheorem{theorem}{Theorem}[section]
\newtheorem{proposition}[theorem]{Proposition}

\newtheorem{corollary}[theorem]{Corollary}
\theoremstyle{definition}
\newtheorem{remark}[theorem]{Remark}
\numberwithin{equation}{section}
\title{Long-root $\SL_2$ dual pairs in split exceptional groups\\
\large Compact types, cone geometry, and theta correspondence}
\author{Gordan Savin}
\date{}
\begin{document}
\maketitle
\begin{abstract}
We study the restriction of the spherical minimal representation $\Pi$
of split real $E_6,E_7,E_8$ to the dual pair $G_1\times G_2$, where
$G_1$ is a long-root $\SL_2(\R)$ and $G_2$ is its centralizer.
We prove that the resulting theta correspondence of representations
is one-to-one. The main ingredients of the proof are lifts of types
of maximal compact subgroups, determined by algebraic-geometric
calculations. Write $\fh=\fk\oplus\fp$ and
$\fg_2=\fk_2\oplus\fp_2$ for the complexified Cartan decompositions
of the ambient group and its centralizer, with maximal compact
subgroups $K$ and $K_2$, respectively. The minimal representation has
a filtration whose associated graded module $R=\gr\Pi$ is the
coordinate ring of the $K$-highest-weight vector cone in $\fp^*$.
Finding generators for the lifts of $K_1=\SO(2)$-types reduces to
computing the affine scheme with coordinate ring $R/(\fp_2)$.
We determine this quotient as a $K_1\times K_2$-module.
The reverse lifts of the $K_2$-types occurring in $R/(\fp_2)$ are
universal $(\mathfrak{sl}_2,\SO(2))$-modules. These calculations lead
to a uniform representation-theoretic description of the three theta
correspondences.
\end{abstract}
\begingroup
\renewcommand{\thefootnote}{}
\endgroup

\section{Introduction}\label{sec:intro}
Let $H$ be the group of real points of a split simply connected algebraic
group of type $E_n$, with $n=6,7,8$, and let $\Pi$ be its spherical minimal
Harish--Chandra module. A long root determines a subgroup
$G_1\simeq\SL_2(\R)$ whose centralizer $G_2$ is, respectively,
$\SL_6(\R)$, split $\Spin_{12}$, or split simply connected $E_7$.
The multiplication map $G_1\times G_2\to H$ has central kernel
$\mu_2$, the group of square roots of $1$. All restrictions below
are understood through this map. Choose compatible maximal compact subgroups
$K\subset H$ and $K_i\subset G_i$, with $K_1=\SO(2)$, and write
$\fh=\fk\oplus\fp$ and $\fg_i=\fk_i\oplus\fp_i$ for their
complexified Cartan decompositions.

For a $(\mathfrak{sl}_2,K_1)$ Harish--Chandra module $\pi$, let
$\pi^\vee$ denote its $K_1$-finite contragredient. We define its big
theta lift by coinvariants:
\[
 \Theta(\pi)=(\Pi\otimes\pi^\vee)_{(\mathfrak{sl}_2,K_1)},
\]
where the subscript denotes coinvariants for the diagonal action.
This construction does not require $\pi$ to be irreducible and defines
a contravariant functor to $(\fg_2,K_2)$-modules. If $\pi$ is
irreducible, then $\pi\boxtimes\Theta(\pi)$ is the largest
$\pi$-isotypic quotient of $\Pi$, recovering the usual definition;
see Baki\'c--Loke--Savin~\cite{BLS25}.

Our purpose is to prove that, for every irreducible $\pi$, the big
lift $\Theta(\pi)$ is a nonzero finite-length $(\fg_2,K_2)$-module
with a unique irreducible quotient, using the ping-pong method
of~\cite{Savin}. To that end we compute the lifts of types of the
maximal compact subgroups of $G_1$ and $G_2$.

Let $v_\omega\in\fp^*$ be a nonzero highest-weight vector for $K$.
The highest-weight vector cone is
\[
 \widehat X=\overline{K_{\C}\cdot v_\omega}\subset\fp^*,
\]
where $K_{\C}$ is the complexification of $K$. Throughout, we identify
$\fp^*$ with $\fp$ using the Killing form. Write
$X=\mathbb P(\widehat X)$ for the projectivization of the cone. With its standard filtration, $\Pi$ has
associated graded module
\[
 R=\gr\Pi=\C[\widehat X].
\]
The quotient $R/(\fp_2)$ is the coordinate ring of the
scheme-theoretic linear section
$\widehat X\cap\fp_2^\perp$. Computing this quotient is the
algebraic-geometric problem that underlies the lifting of types.

Two calculations drive the argument. First, $R/(\fp_2)$ describes finite
generating sets for the lifts of the characters of $\SO(2)$. More precisely,
as an $\SO(2)\times K_2$-module this quotient decomposes as
\[
 R/(\fp_2)=\bigoplus_{k\in\Z}\chi_k\boxtimes Q_k,
\]
where each $Q_k$ is a finite-dimensional graded $K_2$-module. Let
$N(k)=\Hom_{K_1}(\chi_k,R)$ be the lift of $\chi_k$ in the graded model.
Since $\fp_2$ has circle weight zero, $Q_k=N(k)/\fp_2N(k)$.
The graded Nakayama lemma shows that homogeneous $K_2$-equivariant
lifts of $Q_k$ to $N(k)$ generate $N(k)$ over $\Sym(\fp_2)$.
We then lift these generators to the corresponding filtration degrees
in $\Hom_{K_1}(\chi_k,\Pi)$. Induction on the filtration shows that
they generate this lift over $U(\fg_2)$. Thus, in the filtered setting,
the lift of $\chi_k$ is a quotient of the universal module induced from
$Q_k$:
\[
 U(\fg_2)\otimes_{U(\fk_2)}Q_k
   \twoheadrightarrow\Hom_{K_1}(\chi_k,\Pi).
\]
The $K_2$-types occurring in these finite-dimensional spaces $Q_k$ are
the \emph{relevant types}. As described in Section~\ref{sec:forward},
they form two ordered sequences $\sigma_0=1,\sigma_{\pm1},
\sigma_{\pm2},\ldots$, with
\[
 Q_{\pm N}=\bigoplus_{r=0}^{\lfloor N/2\rfloor}
                  \sigma_{\pm(N-2r)}\qquad(N\geq0).
\]
Ordering each sequence by its nonnegative index, $Q_k$ has $\sigma_k$
as its largest $K_2$-type.

The second calculation describes reverse lifts of the relevant types
$\sigma_k$. We have
\[
 R^{K_2}=\C[a,b]=\Sym(\fp_1),\qquad
 \fp_1=\C a\oplus\C b,\qquad \wt(a)=-2,\quad\wt(b)=2.
\]
Here $\wt(a)$ and $\wt(b)$ denote the weights for the circle
$K_1=\SO(2)$. Since $R$ is a domain, it is a torsion-free
$R^{K_2}$-module.
The multiplicity space $M_k=\Hom_{K_2}(\sigma_k,R)$ is therefore
torsion-free as well. The calculation in Section~\ref{sec:reverse} gives
its Poincar\'e series:
\[
 P_{M_k}(t,z)
   =\frac{t^{|k|}z^k}{(1-tz^{-2})(1-tz^2)},
\]
where $t$ records degree and $z$ records circle weight. A nonzero element
of $M_k$ in degree $|k|$ and weight $k$ generates a free rank-one
$R^{K_2}$-submodule, by torsion-freeness. The displayed series shows
that this submodule exhausts $M_k$, so $M_k$ is itself free of rank one
over $R^{K_2}$.
Passing back to $\Pi$, the reverse lift of $\sigma_k$ is the universal
module induced from the $\SO(2)$-type $\chi_k$:
\[
 \Hom_{K_2}(\sigma_k,\Pi)\simeq
 U(\mathfrak{sl}_2)\otimes_{U(\mathfrak{so}_2)}\C_k.
\]

With this indexing, $\sigma_k$ has its primitive occurrence in degree
$|k|$ and circle weight $k$. The representations $\sigma_k$ are listed in
Section~\ref{sec:forward}. Our main conclusion is the following.

\begin{theorem}\label{thm:main}
For each of the three dual pairs and every irreducible
$(\mathfrak{sl}_2,\SO(2))$-module $\pi$, its big theta lift $\Theta(\pi)$
is nonzero and of finite length, with a unique irreducible quotient
$\theta(\pi)$. If $k$ is a weight of $\pi$ of minimal absolute value,
then $\Theta(\pi)$ is generated by one copy of $\sigma_k$. The type
$\sigma_k$ occurs in $\theta(\pi)$ with multiplicity one, and no relevant
type $\sigma_l$ with $|l|<|k|$ occurs there.

The reverse big lift $\Theta(\theta(\pi))$ has $\pi$ as its unique
irreducible quotient. Consequently $\pi\mapsto\theta(\pi)$ is a
bijection from $\Irr(\mathfrak{sl}_2,\SO(2))$ onto the irreducible
$G_2$-modules that occur as partners in this theta correspondence, and
\[
 \dim\Hom_{G_1\times G_2}(\Pi,\pi\boxtimes\theta(\pi))=1.
\]
Here and below homomorphisms intertwine the indicated Lie algebra
and compact-group actions.
\end{theorem}

The theorem concerns the irreducible quotient of the big lift; it does
not assert that the big lift itself is irreducible. A minimal weight is
allowed to have either sign. If both $1$ and $-1$ are minimal weights of
$\pi$, either can be used: its small lift contains both $\sigma_1$ and
$\sigma_{-1}$ with multiplicity one and is generated by either.

We use the branching rules of Loke--Savin~\cite[Proposition~5.1]{LS19}
to prove the invariant and multiplicity identities needed for the present pairs. The minimal
representation and its compact spectrum are supplied by
Brylinski--Kostant~\cite{BK}; see also~\cite{GS}.
For related exceptional dual pairs and their rational forms, see~\cite{LS15}.

Section~\ref{sec:setup} fixes the groups and cone models.
Sections~\ref{sec:forward} and~\ref{sec:reverse} give the forward and
reverse compact-type lifts. Section~\ref{sec:pingpong} proves the theorem.
Sections~\ref{sec:e6}--\ref{sec:e8} establish the three linear-section
calculations, including their scheme structures. In each case the projective
scheme of $R/(\fp_2)$ has two disjoint components, each a projective cone.
The bases of these cones vary with the exceptional group.
They are the projectivizations of the affine cones of symmetric $4\times4$
matrices of rank at most one, arbitrary $4\times4$ matrices of rank at most
one, and skew-symmetric $8\times8$ matrices of rank at most two, respectively.
These give the quadratic Veronese variety $\nu_2(\PP^3)$, the Segre variety
$\PP^3\times\PP^3$, and the Grassmannian $\Gr(2,8)$ in its Pl\"ucker embedding.
The three matrix spaces realize the simple complex Jordan algebras of rank
four occurring in the Brylinski--Kostant construction~\cite{BK}: the
complexifications of the Jordan algebras of $4\times4$ Hermitian matrices
over $\R$, $\C$, and the quaternions, respectively.
In each case the affine cone is the locus of elements of Jordan rank at
most one; in the skew-symmetric realization, Jordan rank one corresponds
to matrix rank two.

A $p$-adic variant of the one-to-one correspondence in
Theorem~\ref{thm:main} is established by Hanzer--Savin~\cite{HS25}.
For the same long-root dual pairs in split $E_6,E_7,E_8$ over a $p$-adic
field, they determine the exceptional theta correspondence and its reverse
lifts; see in particular Sections~5 and~6 of that paper.

\paragraph{Acknowledgments.}
I thank Wee Teck Gan for reading a draft of this paper. I gratefully
acknowledge support from the Simons Foundation through gift No.~946504,
and thank the University of Utah for providing access to OpenAI. 
The paper was written with the help of ChatGPT-6 Astra. 
In addition to harmless assistance with writing, the software was used to compute the ring 
$R/(\fp_2)$. The write up of this result is self contained: 
calculating $R/(\fp_2)$ involves restricting the Pl\"ucker relations to $\fp_2^{\perp}$, 
a pleasant exercise easily done by hand, in this particular case.

\section{The exceptional groups, dual pairs, and graded models}\label{sec:setup}
We describe the three dual pairs, their Cartan decompositions, and the
graded models of the minimal representations.

\subsection{Cartan decompositions and compact groups}
All Lie algebras and algebraic varieties in the calculations are over
$\C$. The real groups are split. We use $\Sp(4)$ for the compact
symplectic group of rank four, with standard complex representation of
dimension eight. The ambient Cartan decompositions are as follows.
\begin{center}
\begin{tabular}{c c c c}
\toprule
$H$ & $K$ & $\fk$ & $\fp$ \\
\midrule
$E_6$ & $\Sp(4)/\mu_2$ & $\mathfrak{sp}_8$ & $\Lambda^4_0\C^8$ \\
$E_7$ & $\SU(8)/\mu_2$ & $\mathfrak{sl}_8$ & $\Lambda^4\C^8$ \\
$E_8$ & $\Spin(16)/\mu_2$ & $\mathfrak{so}_{16}$ & $S^+$ \\
\bottomrule
\end{tabular}
\end{center}
The subscript $0$ means the kernel of symplectic contraction. In the last
row $\mu_2$ is the kernel of the chosen half-spin representation $S^+$;
the compact group is the corresponding semispin quotient.
The dimensions of $\fp$ are $42,70,128$.

For $G_1\simeq\SL_2(\R)$, fix its maximal compact subgroup
$K_1=G_1\cap K\simeq\SO(2)$. Choose $H_0\in\fk_1$ so that
$iH_0$ is the real rotation generator, with period $2\pi$, and put
\[
 k_\theta:=\exp(i\theta H_0)\in K_1,\qquad
 \chi_m(k_\theta)=e^{im\theta},\quad m\in\Z.
\]
Thus $H_0$ acts on the character $\chi_m$ by $m$.
The complexified Cartan decomposition of $\fg_1$ is
\begin{equation}\label{eq:p1}
 \fg_1=\fk_1\oplus\fp_1,\qquad
 \fk_1=\C H_0,\qquad \fp_1=\C E\oplus\C F,
 \quad [H_0,E]=2E,\quad[H_0,F]=-2F.
\end{equation}
Thus $\fp_1=\chi_2\oplus\chi_{-2}$ with the usual integer weights of
$\SO(2)\subset\SL_2(\R)$. For the centralizers we have
\begin{center}
\begin{tabular}{c c c c}
\toprule
$H$ & $G_2$ & $K_2$ & $\fp_2$ \\
\midrule
$E_6$ & $\SL_6(\R)$ & $\SO(6)\simeq\SU(4)/\mu_2$
 & $\Sym^2_0\C^6$ \\
$E_7$ & split $\Spin_{12}$ & $(\SU(4)\times\SU(4))/\mu_2^{\rm diag}$
 & $\Lambda^2\C^4\otimes\Lambda^2\C^4$ \\
$E_8$ & split $E_7$ & $\SU(8)/\mu_2$ & $\Lambda^4\C^8$ \\
\bottomrule
\end{tabular}
\end{center}
The corresponding $\fk_2$ are $\mathfrak{so}_6$,
$\mathfrak{so}_6\oplus\mathfrak{so}_6$, and $\mathfrak{sl}_8$.
The dimensions of $\fp_2$ are $20,36,70$.
The adjoint branching under $\fg_1\oplus\fg_2$ is
\begin{align}
 \mathfrak e_6&=(\mathfrak{sl}_2,1)\oplus(1,\mathfrak{sl}_6)
                  \oplus(\C^2\otimes\Lambda^3\C^6),\notag\\
 \mathfrak e_7&=(\mathfrak{sl}_2,1)\oplus(1,\mathfrak{so}_{12})
                  \oplus(\C^2\otimes S_{12}),\label{eq:adjoint}\\
 \mathfrak e_8&=(\mathfrak{sl}_2,1)\oplus(1,\mathfrak e_7)
                  \oplus(\C^2\otimes V_{56}).\notag
\end{align}
Here $S_{12}$ is a half-spin module of dimension $32$, and $V_{56}$ is
the $56$-dimensional fundamental $E_7$-module. These descriptions identify
the long-root centralizers directly from the adjoint branching.

Under the inclusion $K_1\subset K$, the element $k_\theta$ is
$\exp_K(i\theta H_0)$, with $H_0$ viewed in $\fk$.
For example, in $E_7$ write $U=V\oplus W$ with
$\dim V=\dim W=4$. The images in $\fk=\mathfrak{sl}(U)$ and
$K=\SU(8)/\mu_2$ are
\begin{align*}
 H_0&=\tfrac12\operatorname{diag}(I_V,-I_W),\\
 k_\theta=\exp_K(i\theta H_0)
   &=[\operatorname{diag}(e^{i\theta/2}I_V,e^{-i\theta/2}I_W)].
\end{align*}
When $\theta$ increases by $2\pi$, the representative in $\SU(8)$
changes by $-I$, which is trivial in $K$. This realizes the period
$2\pi$ of $K_1$. The auxiliary spaces $V,W$ have $H_0$-weights
$1/2,-1/2$ on the cover, while the actual $K_1$-modules have integer
weights. In particular the weights of $\fp_1$ are $2,-2$, as in
\eqref{eq:p1}.

\subsection{The minimal representation and its filtration}
Let $v_0$ be a nonzero spherical vector in $\Pi$, and let $\omega$ be
the highest weight of $\fp$ as a $K$-module. The compact spectrum is
$\Pi|_K=\bigoplus_{d\geq0}V_K(d\omega)$, with multiplicity one.
Its standard filtration has
\[
 F_d\Pi=\bigoplus_{0\leq e\leq d}V_K(e\omega),\qquad
 \gr\Pi=R=\bigoplus_{d\geq0}H^0(X,\OO_X(d)).
\]
The associated graded action of $\Sym(\fp)$ is multiplication on the cone ring,
using the Killing-form identification of $\fp$ with $\fp^*$.
These are the graded models of the spherical minimal representations
in~\cite{BK}. Explicitly:
\begin{center}
\begin{tabular}{c c c c}
\toprule
$H$ & $X$ & $R_d$ & $\dim X$ \\
\midrule
$E_6$ & $\LG(4,8)$ & $V_{C_4}(d\omega_4)^*$ & $10$ \\
$E_7$ & $\Gr(4,8)$ & $\bS_{(d^4)}(\C^8)^*$ & $16$ \\
$E_8$ & $\OGr^+(8,16)$ & $V_{D_8}(d\omega_8)^*$ & $28$ \\
\bottomrule
\end{tabular}
\end{center}
The first two embeddings are Pl\"ucker embeddings. For $E_8$ we label
the selected half-spin weight by $\omega_8$ and use the spinor embedding.
In each case $R$ is generated in degree one and is an integral domain.

\subsection{The circle decomposition and the Jordan algebra}
Let $J$ be the rank-four complex Jordan algebra associated with $H$
in the introduction. Write $\C_m$ for a one-dimensional space on which
$K_1$ acts by $\chi_m$, and write $J_1$ and $J^*_{-1}$ for $J$ and
its dual with circle weights $1$ and $-1$, respectively.
The circle decomposition of $\fp$ is
\begin{equation}\label{eq:degreeone}
 \fp=\C_{-2}\oplus J^*_{-1}\oplus\fp_2\oplus J_1\oplus\C_2.
\end{equation}
The middle term has circle weight zero, and the two $K_2$-trivial
scalar lines form $\fp_1$. The Killing form pairs opposite weight
spaces and identifies $\fp_2$ with linear functions on the cone;
$(\fp_2)$ denotes the ideal they generate in $R$.

\section{Lifts of the \texorpdfstring{$\SO(2)$}{SO(2)}-types: the common pattern}\label{sec:forward}
We determine generators for the lifts of $\SO(2)$-types in $R$ and
$\Pi$ using the quotient $R/(\fp_2)$.

\subsection{Cone bases and relevant types}
Put $Q=R/(\fp_2)$. We call a $K_2$-type \emph{relevant} if it occurs
in $Q$. By \eqref{eq:degreeone}, setting $\fp_2$ equal to zero cuts
the highest-weight vector cone by
$\fp_2^\perp=(\C_2\oplus J_1)\oplus(\C_{-2}\oplus J^*_{-1})$.
Sections~\ref{sec:e6}--\ref{sec:e8} show that this section is reduced
and that $\operatorname{Proj}Q=X_+\sqcup X_-$, where
\[
 \begin{aligned}
 X_+&=\{[\lambda:u]\in\PP(\C_2\oplus J_1):
                         \operatorname{rank}_J u\leq1\},\\
 X_-&=\{[\lambda:v]\in\PP(\C_{-2}\oplus J^*_{-1}):
                         \operatorname{rank}_{J^*}v\leq1\}.
 \end{aligned}
\]
Here the ranks are Jordan ranks in the two dual matrix realizations.
The vertex of $X_+$ is $[1:0]$, representing the line $\C_2$.
Its base is the rank-one locus in $\PP(J_1)=\PP(J)$, embedded by
$[u]\mapsto[0:u]$. The cone $X_+$ is the union of the projective lines
$[\lambda:\mu u]$, with $[\lambda:\mu]\in\PP^1$, joining the vertex
to these base points. Similarly, $X_-$ has vertex the line $\C_{-2}$
and base the rank-one locus in $\PP(J^*)$. The vertices belong to the
two components, while the affine cones meet only at the origin.
The bases, with their projective embeddings, are as follows.
\begin{center}
\begin{tabular}{c c c c}
\toprule
$H$ & base of $X_+$ & base of $X_-$ & $\dim X_\pm$ \\
\midrule
$E_6$ & $\nu_2(\PP(V^*))$ & $\nu_2(\PP(V))$ & $4$ \\
$E_7$ & $\PP(V^*)\times\PP(W)$ & $\PP(V)\times\PP(W^*)$ & $7$ \\
$E_8$ & $\Gr(2,V^*)$ & $\Gr(2,V)$ & $13$ \\
\bottomrule
\end{tabular}
\end{center}
Here $\dim V=4$ in $E_6$, $\dim V=\dim W=4$ in $E_7$, and
$\dim V=8$ in $E_8$. The embeddings are the second Veronese, the Segre,
and the Pl\"ucker embedding, respectively. Accordingly, the $K_2$-module
$J$ is $\Sym^2V^*$, $V^*\otimes W$, or $\Lambda^2V^*$, respectively.

Let $\mathcal S^+=\bigoplus_{j\geq0}\mathcal S_j^+$ and
$\mathcal S^-=\bigoplus_{j\geq0}\mathcal S_j^-$ be the coordinate
rings of the rank-at-most-one cones in $J$ and $J^*$, respectively,
that is, the affine cones over the two bases in the table.
By the standard homogeneous coordinate-ring decompositions
\cite{Fulton,FH}, each degree-$j$ piece is irreducible under $K_2$;
denote its $K_2$-type by $\tau_j^\pm$. Including the circle action, we have
\[
 \mathcal S_j^+\simeq\chi_{-j}\boxtimes\tau_j^+,
 \qquad
 \mathcal S_j^-\simeq\chi_j\boxtimes\tau_j^-,
 \qquad \tau_0^\pm=1.
\]
The coordinate weights are opposite to the vector weights on
$J_1$ and $J^*_{-1}$. Besides the common trivial type, these are the
relevant types; for $j\geq1$ they are listed below.
\begin{center}
\begin{tabular}{c c c}
\toprule
$H$ & $\tau_j^+$ & $\tau_j^-$ \\
\midrule
$E_6$ & $\Sym^{2j}V$ & $\Sym^{2j}V^*$ \\
$E_7$ & $\Sym^jV\otimes\Sym^jW^*$ & $\Sym^jV^*\otimes\Sym^jW$ \\
$E_8$ & $\bS_{(j,j)}V=V_{j\omega_2}$ & $(\bS_{(j,j)}V)^*=V_{j\omega_6}$ \\
\bottomrule
\end{tabular}
\end{center}
The notation $\tau_j^\pm$ records only the $K_2$-type; the circle
weights are recorded separately above. All these types descend through
the compact quotients in Section~\ref{sec:setup}.
Geometrically, they are the irreducible constituents of
$H^0(X_\pm,\OO_{X_\pm}(d))$, for $d\geq0$, where $\OO(1)$ is the
ambient Pl\"ucker line bundle in $E_6,E_7$ and the spinor line bundle
in $E_8$. Sections~\ref{sec:e6}--\ref{sec:e8} establish the ring
descriptions. The section interpretation follows from projective
normality of the Veronese, Segre, and Grassmannian embeddings and
their projective cones.

In the formulas below, $a,b\in R^{K_2}$ are coordinates on the scalar
summands $\C_2,\C_{-2}$, extended by zero on the other summands;
we use the same notation for their images in $Q$. They have degree
one and weights $-2,2$, since coordinate weights are opposite to vector
weights.

\begin{theorem}\label{thm:quotient}
In each case $Q=R/(\fp_2)$, the coordinate ring of the linear section,
is reduced and has graded character
\begin{equation}\label{eq:quotientchar}
 \ch_{t,z}Q=
 \frac{\sum_{j\geq0}t^jz^{-j}\ch(\tau_j^+)}{1-tz^{-2}}
 +\frac{\sum_{j\geq0}t^jz^j\ch(\tau_j^-)}{1-tz^2}-1.
\end{equation}
More precisely, as a graded $K_1\times K_2$-algebra,
\begin{equation}\label{eq:quotientmodule}
 Q\simeq\mathcal S^+[a]\times_\C\mathcal S^-[b].
\end{equation}
The fiber product consists of pairs with the same constant term.
In particular $Q$ is admissible even under $\SO(2)$ alone. For $N>0$,
\begin{equation}\label{eq:forwardgenerators}
 Q_{-N}=\bigoplus_{r=0}^{\lfloor N/2\rfloor}\tau_{N-2r}^+,
 \qquad Q_N=\bigoplus_{r=0}^{\lfloor N/2\rfloor}\tau_{N-2r}^-,
 \qquad Q_0=1.
\end{equation}
The summand indexed by $r$ has degree $N-r$.
\end{theorem}

The scheme and graded-ring calculations proving this theorem are given in
Sections~\ref{sec:e6}--\ref{sec:e8}. To make the circle action entirely
explicit, forgetting the degree, \eqref{eq:quotientmodule} means
\begin{equation}\label{eq:explicitchars}
 Q=\bigoplus_{q\in\Z}\chi_{2q}\boxtimes1
 \ \oplus\!\bigoplus_{j\geq1,\,q\geq0}\chi_{-j-2q}\boxtimes\tau_j^+
 \ \oplus\!\bigoplus_{j\geq1,\,q\geq0}\chi_{j+2q}\boxtimes\tau_j^-.
\end{equation}
The degree in the first sum is $|q|$, and in the other sums it is $j+q$.
In particular the weight-zero lift in this quotient is just $1$.

\subsection{Passage to the minimal representation}
For $k\in\Z$ set
\[
 N(k)=\Hom_{K_1}(\chi_k,R),\qquad
 A(k)=\Hom_{K_1}(\chi_k,\Pi).
\]
The graded $\Sym(\fp_2)$-module $N(k)$ retains the grading of $R$:
its degree-$d$ part is $\Hom_{K_1}(\chi_k,R_d)$. The filtered module
$A(k)$ is a $(\fg_2,K_2)$-module.
Define the signed notation
\begin{equation}\label{eq:sigma}
 \sigma_0=1,\qquad \sigma_{-j}=\tau_j^+,\qquad
 \sigma_j=\tau_j^-\quad(j\geq1).
\end{equation}
Thus the index of $\sigma_k$ is its primitive circle weight.

For $k\ne0$, put
\begin{equation}\label{eq:Fk}
 F(k)=\bigoplus_{r=0}^{\lfloor |k|/2\rfloor}
        \sigma_{\operatorname{sgn}(k)(|k|-2r)},\qquad F(0)=1.
\end{equation}
When $|k|$ is even the last summand is $\sigma_0$, included only once.

\begin{corollary}\label{cor:forward}
The module $A(k)$ is generated over $U(\fg_2)$ by filtered lifts of
$F(k)$. The summand indexed by $r$ lifts in filtration degree $|k|-r$.
In particular every $A(k)$ is finitely generated, and $A(0)$ is generated
by the spherical vector.
\end{corollary}
\begin{proof}
Taking a fixed compact type commutes with associated graded because
compact representations are semisimple. Thus $\gr A(k)=N(k)$. Since
$\fp_2$ has circle weight zero,
$N(k)/\fp_2N(k)=Q_k=F(k)$ with the degrees specified above.
Choose homogeneous lifts of a basis of this quotient. Induction on
degree shows that they generate $N(k)$ over $\Sym(\fp_2)$: after
subtracting their constant coefficients, an element lies in $\fp_2N(k)$,
and its coefficients have smaller degree. Lift these generators
$K_2$-equivariantly to $A(k)$ and apply the same filtration induction
using $U(\fg_2)$. This proves the claim.
\end{proof}

\section{Reverse lifts of the relevant types}\label{sec:reverse}
For a $K_2$-type $\tau$ define
\[
 M(\tau)=\Hom_{K_2}(\tau,R),\qquad
 B(\tau)=\Hom_{K_2}(\tau,\Pi).
\]
Their gradings and filtrations retain the circle action. We compute only
the reverse lifts needed for the relevant types.

\begin{theorem}\label{thm:reverse}
In each of the three cases,
\begin{equation}\label{eq:invariants}
 R^{K_2}=\C[a,b]=\Sym(\fp_1).
\end{equation}
For every signed relevant type $\sigma_k$, its graded multiplicity space
$M(\sigma_k)$ is free of rank one over this invariant ring, with generator in degree
$|k|$ and weight $k$. Its Poincar\'e series is
\begin{equation}\label{eq:reverseP}
 P_{\sigma_k}(t,z)=
 \frac{t^{|k|}z^k}{(1-tz^{-2})(1-tz^2)}.
\end{equation}
Moreover, as a $(\mathfrak{sl}_2,\SO(2))$-module,
\begin{equation}\label{eq:universal}
 B(\sigma_k)\simeq I(k):=
 U(\mathfrak{sl}_2)\otimes_{U(\mathfrak{so}_2)}\C_k,
\end{equation}
with its generating line in filtration degree $|k|$.
\end{theorem}

\subsection{Branching and the Poincar\'e series}
We use the restriction of $R_d$ to $K_1\times K_2$, retaining both
its circle weights and its $K_2$-types. The branching rules are given
by Loke--Savin~\cite[Proposition~5.1]{LS19}, with $m=4$ in each
of its three cases. We state only their consequence for the relevant
types, in the notation and integer circle-weight normalization of
Section~\ref{sec:forward}.

For every $k\in\Z$ and $d\geq0$, the $K_1$-module of multiplicities is
\begin{equation}\label{eq:relevantbranch}
 \Hom_{K_2}(\sigma_k,R_d)
   =\bigoplus_{\substack{r,s\geq0\\r+s=d-|k|}}
       \chi_{k-2r+2s}.
\end{equation}
The sum is empty when $d<|k|$. Thus $\sigma_k$ first occurs in degree
$|k|$, with circle weight $k$. In degree $|k|+n$, its occurrences
have the distinct weights
\[
 k-2n,\ k-2n+4,\ \ldots,\ k+2n,
\]
each with multiplicity one. Equivalently, the
$\sigma_k$-isotypic part of $R_d$, as a $K_1\times K_2$-module,
is the sum of $\chi_{k-2r+2s}\boxtimes\sigma_k$ over the same indices.
This statement includes $k=0$, for which $\sigma_0=1$.

Writing $d=|k|+r+s$, we obtain at once
\[
 P_{\sigma_k}(t,z)
 =\sum_{r,s\geq0}t^{|k|+r+s}z^{k-2r+2s}
 =\frac{t^{|k|}z^k}{(1-tz^{-2})(1-tz^2)}.
\]
This proves \eqref{eq:reverseP} uniformly in the three cases.

\subsection{From the series to free and universal modules}
Taking $k=0$ in \eqref{eq:relevantbranch} gives
\[
 \ch_{t,z}R^{K_2}=\frac1{(1-tz^{-2})(1-tz^2)}.
\]
The elements $a,b$ are nonzero in the domain $R$. Their monomials are
nonzero and have distinct pairs of degree and weight. They are therefore
linearly independent and realize every occurrence in this series.
This proves \eqref{eq:invariants}, including the algebraic independence
of $a,b$.

Choose a nonzero homogeneous equivariant map
$g_k:\sigma_k\to R_{|k|}$ of circle weight $k$.
Multiplication defines
\[
 \C[a,b]\,g_k\longrightarrow M(\sigma_k).
\]
It is injective: a nonzero polynomial cannot annihilate the nonzero
image of $g_k$ in the domain $R$. The source and target have the same
degree-and-weight series \eqref{eq:reverseP}, so the map is surjective.
This proves the rank-one assertion. A Poincar\'e series by itself would
not establish freeness; the injective multiplication map is essential.

Now lift $g_k$ $K_1$-equivariantly to $B(\sigma_k)$.
Filtered generation gives a surjection $I(k)\to B(\sigma_k)$.
PBW gives
\[
 \gr I(k)=\Sym(\fp_1)\otimes\C_k,
\]
with the generating line placed in degree $|k|$.
The associated graded map is the rank-one isomorphism just proved, so
the filtered map is an isomorphism. This finishes the proof of
Theorem~\ref{thm:reverse}.

Write $\pi_k=\Hom_{K_1}(\chi_k,\pi)$
for the weight-$k$ multiplicity space. The universal property we need is
\begin{equation}\label{eq:universalhom}
 \Hom_{\fg_1,K_1}(I(k),\pi)=\pi_k.
\end{equation}

\begin{remark}\label{rem:allreverse}
The same argument applies to every $K_2$-type $\tau$ occurring in $R$.
Indeed, fixing its highest-weight differences in
\cite[Proposition~5.1]{LS19} leaves two nonnegative parameters. There
is a unique first occurrence, in degree $d_\tau$ and circle weight
$w_\tau$, and
\[
 \Hom_{K_2}(\tau,R_d)
 =\bigoplus_{\substack{r,s\geq0\\r+s=d-d_\tau}}
       \chi_{w_\tau-2r+2s},\qquad
 P_\tau(t,z)=\frac{t^{d_\tau}z^{w_\tau}}
 {(1-tz^{-2})(1-tz^2)}.
\]
Since $R$ is a domain, multiplication by a nonzero first-occurrence
map identifies $M(\tau)$ with a free rank-one $R^{K_2}$-module.
The filtered argument then gives
$B(\tau)\simeq I(w_\tau)$, with its generating line in degree
$d_\tau$. For $\tau=\sigma_k$ these parameters are $d_\tau=|k|$
and $w_\tau=k$; in general $d_\tau$ need not equal $|w_\tau|$.
\end{remark}

\section{Ping-pong and the uniform theta correspondence}\label{sec:pingpong}
We prove finite length of the big theta lifts and use the compact-type
calculations to establish Theorem~\ref{thm:main}.

\subsection{Finiteness and compact-type reciprocity}
For an irreducible $(\fg_1,K_1)$-module $\pi$, recall that
$\pi\boxtimes\Theta(\pi)$ is the maximal $\pi$-isotypic quotient
of $\Pi$. Define the reverse big lift
of an irreducible $(\fg_2,K_2)$-module $\rho$ similarly, with quotient
$\Theta(\rho)\boxtimes\rho$. Compact-type reciprocity reads
\begin{align}
 \Hom_{K_2}(\Theta(\pi),\tau)
   &\simeq\Hom_{(\fg_1,K_1)\times K_2}(\Pi,\pi\boxtimes\tau)
    \simeq\Hom_{\fg_1,K_1}(B(\tau),\pi),\label{eq:reciprocity1}\\
 \Hom_{K_1}(\Theta(\rho),\chi_l)
   &\simeq\Hom_{K_1\times(\fg_2,K_2)}(\Pi,\chi_l\boxtimes\rho)
    \simeq\Hom_{\fg_2,K_2}(A(l),\rho).\label{eq:reciprocity2}
\end{align}
The first isomorphisms use the maximal isotypic quotients, and the second
use the compact isotypic components.
For a locally finite $K$-module $M$ and an irreducible $K$-type $\eta$,
$\Hom_K(M,\eta)\simeq\Hom_K(\eta,M)^*$, with algebraic dual, so
finite-dimensionality is equivalent to finite multiplicity.

\begin{proposition}\label{prop:finiteness}
Both big theta lifts have finite length whenever the input is
irreducible. Moreover, for every relevant type,
\begin{equation}\label{eq:exactmultiplicity}
 \dim\Hom_{K_2}(\sigma_l,\Theta(\pi))=\dim\pi_l.
\end{equation}
In particular $\Theta(\pi)\ne0$ for every irreducible $\pi$.
\end{proposition}
\begin{proof}
For any $K_2$-type $\tau$ occurring in $R$,
Remark~\ref{rem:allreverse} gives $B(\tau)\simeq I(w_\tau)$.
Thus \eqref{eq:reciprocity1} and \eqref{eq:universalhom} identify
$\Hom_{K_2}(\Theta(\pi),\tau)$ with $\pi_{w_\tau}$, which is
finite-dimensional. Types absent from $R$ have zero multiplicity.
Likewise, Corollary~\ref{cor:forward} and \eqref{eq:reciprocity2} give
\[
 \dim\Hom_{K_1}(\Theta(\rho),\chi_l)
 =\dim\Hom_{\fg_2,K_2}(A(l),\rho)
 \leq\dim\Hom_{K_2}(F(l),\rho)<\infty.
\]
The preceding duality now proves that both big lifts are admissible.
To obtain finite length, we also verify finite generation using the
compact-type lifts.

Every circle weight space of an irreducible
$(\mathfrak{sl}_2,\SO(2))$-module is at most one-dimensional.
For any occurring $k$, taking weight $k$ in
$\Pi\twoheadrightarrow\pi\boxtimes\Theta(\pi)$ gives
$A(k)\twoheadrightarrow\Theta(\pi)$, proving finite generation.
For the reverse direction, choose a $K_2$-type $\tau$ occurring in
$\rho$. Taking its isotypic multiplicity space gives a surjection
$B(\tau)\twoheadrightarrow
\Theta(\rho)\otimes\Hom_{K_2}(\tau,\rho)$; projection to a nonzero
factor proves finite generation of $\Theta(\rho)$.
Admissible finitely generated modules have finite length
\cite{Wallach}. Finally \eqref{eq:exactmultiplicity} follows from
\eqref{eq:reciprocity1}, Theorem~\ref{thm:reverse}, and
\eqref{eq:universalhom}. Some $\pi_l$ is nonzero, proving nonvanishing.
\end{proof}

\subsection{The two passes}
For an arbitrary pair $(\fg,K)$, let $M$ be a nonzero $(\fg,K)$-module
of finite length as a $U(\fg)$-module. If $M$ is generated over
$U(\fg)$ by one copy of an irreducible $K$-type $\eta$ and that type
has multiplicity one in $M$, then $M$ has a unique irreducible quotient. Indeed every nonzero quotient must contain the
image of its generating copy of $\eta$. The semisimple head therefore
has only one irreducible summand, with multiplicity one.

\begin{proof}[Proof of Theorem~\ref{thm:main}]
Choose a circle weight $k$ of $\pi$ with minimal absolute value.
Equation~\eqref{eq:exactmultiplicity} gives
\begin{equation}\label{eq:minimalmultiplicity}
 [\Theta(\pi):\sigma_k]_{K_2}=1,\qquad
 [\Theta(\pi):\sigma_l]_{K_2}=0\quad (|l|<|k|).
\end{equation}
The surjection $A(k)\twoheadrightarrow\Theta(\pi)$ and
Corollary~\ref{cor:forward} show that the big lift is generated by the
image of $F(k)$. All its summands except $\sigma_k$ vanish by
\eqref{eq:minimalmultiplicity}. Thus $\Theta(\pi)$ is generated by
one copy of $\sigma_k$, and the observation gives its unique irreducible
quotient $\rho=\theta(\pi)$. Its generating type survives and still
has multiplicity one. Since maps to $\rho$ are determined by their
restriction to this generating type, $\dim\Hom_{\fg_2,K_2}(\Theta(\pi),\rho)=1$,
which gives the multiplicity-one assertion of the theorem.
The smaller relevant types are absent.

For the return pass, the quotient
$\Pi\twoheadrightarrow\pi\boxtimes\rho$ makes $\pi$ a quotient of
$\Theta(\rho)$. By \eqref{eq:reciprocity2},
\[
 \dim\Theta(\rho)_k
   =\dim\Hom_{\fg_2,K_2}(A(k),\rho)
   \leq\dim\Hom_{K_2}(F(k),\rho)=1.
\]
The inequality follows by restricting a homomorphism to the chosen
generators. The quotient onto $\pi$ makes the dimension at least one,
so it equals one. On the other hand, taking the $\sigma_k$-multiplicity
space of the reverse maximal quotient gives
\[
 I(k)=B(\sigma_k)\twoheadrightarrow
 \Theta(\rho)\otimes\Hom_{K_2}(\sigma_k,\rho)
 \simeq\Theta(\rho).
\]
Therefore $\Theta(\rho)$ is generated by its single weight-$k$ line.
Every nonzero quotient contains that line, so the same observation
proves that its unique irreducible quotient is $\pi$.

\end{proof}

\section{The linear section in type \texorpdfstring{$E_6$}{E6}}\label{sec:e6}
We compute $Q=R/(\fp_2)$ for the long-root $\SL_2\times\SL_6$
pair in $E_6$. We identify its two Veronese cones and their scheme
structure, then read off the graded $\SO(2)\times K_2$-decomposition
from their coordinate rings to prove Theorem~\ref{thm:quotient}.

\subsection{The two Veronese cones and their equations}
Let $U=V\oplus V^*$, $\dim V=4$, with symplectic form
\[
 \langle(v,\phi),(w,\psi)\rangle=\psi(v)-\phi(w).
\]
A point $L\in\LG(4,U)$ is a four-dimensional Lagrangian subspace of $U$.
The ambient Cartan space is $\fp=\Lambda^4_0U$, the kernel of
symplectic contraction, and it contains the Pl\"ucker vectors of
$\LG(4,U)$. As an $\SU(V)$-module it decomposes as
\[
 \fp\simeq\C\oplus\Sym^2V^*
       \oplus\underbrace{\bS_{(2,2)}V}_{\fp_2}
       \oplus\Sym^2V\oplus\C.
\]
The summands are ordered by the number $4,3,2,1,0$ of $V$ factors
in $\Lambda^4(V\oplus V^*)$. Their dimensions are $1,10,20,10,1$,
and their circle weights are $2,1,0,-1,-2$.
Determinant characters have been suppressed in this display because
they are trivial on $\SU(V)$.

To see the decomposition, the two extreme terms are $\det V$ and
$\det V^*$. The term $\Lambda^3V\otimes V^*$ identifies with
$\det V\otimes(V^*\otimes V^*)$; contraction vanishes on the
symmetric part and is injective on the alternating part. Its kernel
is therefore $\det V\otimes\Sym^2V^*$, and the opposite term is dual.
On the middle term $\Lambda^2V\otimes\Lambda^2V^*$, contraction
maps onto $V\otimes V^*$ and has the irreducible $20$-dimensional
kernel $\bS_{(2,2)}V$. Equivalently, this is
$\Sym^2_0(\Lambda^2V)\simeq\Sym^2_0\C^6$, the Cartan space $\fp_2$
of the centralizer.

The coordinate functions have the opposite circle weights.
Thus $\fp_2=0$ sets the middle Pl\"ucker coordinates to zero,
leaving the four spaces with weights $-2,-1,1,2$.
To describe the resulting four-planes, put
$r=\dim(L\cap V)$ and $s=\dim(L\cap V^*)$.
Complete bases of these intersections by vectors $v_i+\phi_i$ in $L$.
The classes of the $v_i$ and of the $\phi_i$ are independent in the
respective quotient spaces. Expanding their wedge shows that the
nonzero components have precisely $r,r+1,\ldots,4-s$ factors from $V$.
The middle component therefore vanishes exactly when $r\geq3$ or
$s\geq3$. We obtain the two loci
\[
 X_+=\{L\in\LG(4,U):\dim(L\cap V)\geq3\},\qquad
 X_-=\{L\in\LG(4,U):\dim(L\cap V^*)\geq3\}.
\]
They are disjoint: the two indicated intersections cannot both have
dimension at least three inside a four-dimensional space.

To understand $X_+$, take $L\ne V$ and put $A=L\cap V$.
Then $\dim A=3$, and $L/A$ is a line in $(V/A)\oplus V^*$
with nonzero projection to $V^*$. Isotropy forces this projection
to annihilate $A$: for $a\in A$ and $v+\phi\in L$, we have
$\langle a,v+\phi\rangle=\phi(a)=0$. Since the annihilator
\[
 D=\operatorname{Ann}(A)=\{\phi\in V^*: \phi(A)=0\}
\]
is itself a line, the projected line must be $D$. The choice of $A$
therefore determines $D$; there is no additional choice of a line in
$V^*$. In particular the base will be $\Gr(3,V)=\PP(V^*)$.

Conversely, fix $A\in\Gr(3,V)$, and put $Q_A=V/A$, so
$D=Q_A^*$. The orthogonal complement for the symplectic form is
$A^\perp=V\oplus D$, and
\[
 A^\perp/A=Q_A\oplus Q_A^*
\]
is a two-dimensional symplectic space. Every line in this space is
isotropic, and its inverse image in $A^\perp$ is a Lagrangian
four-plane containing $A$. These four-planes form a $\PP^1$.
The line $Q_A$ gives $V$, and the line $Q_A^*$ gives $A\oplus D$.
Every other member also has intersection exactly $A$ with $V$.

In the Pl\"ucker embedding this $\PP^1$ is the line joining
$V$ to $A\oplus D$. Indeed, if $u_1,u_2,u_3$ is a basis of $A$,
$v\in V\setminus A$, and $0\ne\phi\in D$, its points are
\[
 [u_1\wedge u_2\wedge u_3\wedge(\lambda v+\mu\phi)],
 \qquad[\lambda:\mu]\in\PP^1.
\]
To see the embedding of the base, use
$\Lambda^3V\simeq\det V\otimes V^*$. Under this identification
$\Lambda^3A$ corresponds to $\det V\otimes D$. The Pl\"ucker
line of $A\oplus D$ is consequently
$\det V\otimes D^{\otimes2}$. As $D$ varies through the lines
in $V^*$, this is exactly the second Veronese embedding
$\nu_2:\PP(V^*)\to\PP(\Sym^2V^*)$, with the harmless common
factor $\det V$. Thus $X_+$ is the projective cone over this
Veronese variety, with vertex $V$. Interchanging $V$ and $V^*$
gives $X_-$ and its vertex $V^*$.

We now check this description scheme-theoretically. Start with the
Pl\"ucker ideal of $\Gr(4,V\oplus V^*)$, generated by its degree-two
part $I_2$. Passing to $\LG(4,U)$ adds the linear contraction equations
just described. On the surviving coordinates they make the two
$4\times4$ arrays $x,y$ symmetric, with the usual cofactor signs.
Killing the middle summand $\fp_2$ amounts to restricting the equations
that kill the full middle coordinate space in the ordinary Grassmannian.

Thus the surviving coordinate spaces are
\[
 \C a,\qquad \mathcal X=\operatorname{span}\{x_{ij}:i\leq j\}
       \simeq\Sym^2V,\qquad
 \mathcal Y=\operatorname{span}\{y_{ij}:i\leq j\}
       \simeq\Sym^2V^*,\qquad \C b,
\]
as $\SU(V)$-modules, with circle weights $-2,-1,1,2$.
Here $x_{ji}=x_{ij}$ and $y_{ji}=y_{ij}$.
Write $\overline I_2$ for the image of the Pl\"ucker quadratic
relation space after imposing these linear equations. It lies in
$S=\C[a,b,x,y]$, and the defining ideal of $Q$ is $(\overline I_2)$.
Let $\overline I_2(x)$ and $\overline I_2(y)$ denote its subspaces
involving only $x$ or only $y$, respectively. By the calculation in
Section~\ref{sec:e7}, they are spanned by the $2\times2$ minors of
the corresponding symmetric matrices. These are spaces of quadratic
polynomials; parentheses will denote the ideals they generate.

The $\SU(V)$-module structure of $\overline I_2(x)$ follows from
multiplication of quadratic polynomials:
\[
 \overline I_2(x)
 =\ker\bigl(\Sym^2(\Sym^2V)\longrightarrow\Sym^4V\bigr)
 =\bS_{(2,2)}V.
\]
Indeed, the standard Schur decomposition is
$\Sym^2(\Sym^2V)=\Sym^4V\oplus\bS_{(2,2)}V$
\cite[Section~1.2]{Ottaviani}. The kernel is spanned by
$x_{ij}x_{kl}-x_{il}x_{kj}$, which vanish under $x_{ij}\mapsto s_is_j$.
It has dimension $55-35=20$ and highest weight $2\omega_2$.
The dual space $\overline I_2(y)=\bS_{(2,2)}V^*$ is also
isomorphic to this irreducible $\SU(4)$-module.

Specializing the generator decomposition of Section~\ref{sec:e7}
now gives
\[
 \overline I_2=\overline I_2(x)\oplus a\mathcal Y
 \oplus\bigl(\C ab\oplus\mathcal X\mathcal Y\bigr)
 \oplus b\mathcal X\oplus\overline I_2(y).
\]
To display its irreducible constituents, let $V_\lambda$ denote
the $\SU(4)$-module of highest weight $\lambda$, with fundamental
weights $\omega_1,\omega_2,\omega_3$. The tensor product in the
middle is
\[
 \mathcal X\mathcal Y\simeq\Sym^2V\otimes\Sym^2V^*
 =\C\oplus V_{\omega_1+\omega_3}
       \oplus V_{2\omega_1+2\omega_3},
\]
with dimensions $1,15,84$. Grouping by the total $V$-degree of the
original Pl\"ucker indices, the quadratic generators are:
\begin{center}
\begin{tabular}{c c c c}
\toprule
total $V$-degree & generator space & $\SU(4)$-module & dimension \\
\midrule
$6$ & $\overline I_2(x)$ & $V_{2\omega_2}$ & $20$ \\
$5$ & $a\mathcal Y$ & $V_{2\omega_3}$ & $10$ \\
$4$ & $\C ab\oplus\mathcal X\mathcal Y$
 & $\C^{\oplus2}\oplus V_{\omega_1+\omega_3}\oplus V_{2\omega_1+2\omega_3}$ & $101$ \\
$3$ & $b\mathcal X$ & $V_{2\omega_1}$ & $10$ \\
$2$ & $\overline I_2(y)$ & $V_{2\omega_2}$ & $20$ \\
\bottomrule
\end{tabular}
\end{center}
The two copies of the trivial representation in degree $4$ are
$\C ab$ and the invariant contraction in $\mathcal X\mathcal Y$.
All dimensions follow from the displayed modules; their sum is $161$.
The degree-two part of $Q$ consequently has dimension
$\binom{23}{2}-161=92$.

The relations in total $V$-degrees $5,4,3$ are exactly
$ay_{ij},ab,x_{ij}y_{kl},bx_{ij}$. They generate $(a,x)(b,y)$,
so they separate the two coordinate blocks, while the remaining
relations impose rank at most one on the symmetric matrix in each
block. Consequently,
\begin{equation}\label{eq:ring6}
 Q=\bigl(\C[x]/(\overline I_2(x))\bigr)[a]\times_\C
       \bigl(\C[y]/(\overline I_2(y))\bigr)[b].
\end{equation}
Here $\C[x]=\C[x_{ij}:i\leq j]$ and similarly for $y$; each map
to $\C$ evaluates at the origin. The substitution $x_{ij}=s_is_j$
identifies the first matrix ring with
\[
 \C[x]/(\overline I_2(x))
   =\C[s_is_j:1\leq i\leq j\leq4]
   =\bigoplus_{d\geq0}\Sym^{2d}V,
\]
where a quadratic monomial in the $s_i$ has degree one in the matrix
ring. This is the homogeneous coordinate ring of the quadratic
Veronese variety $\nu_2(\PP(V^*))$. Its affine spectrum consists
of symmetric matrices of rank at most one; its projectivization
consists of their nonzero rank-one classes. The second factor is dual.
Adjoining the free vertex coordinates $a,b$ gives the two projective
cones. The fibre product identifies their constant terms: the affine
cones meet only at the origin, while their projectivizations are disjoint.
Both factors are polynomial extensions of reduced Veronese rings,
which proves the scheme assertion.

The degree-$j$ piece of the first matrix ring is therefore
$\mathcal S_j^+\simeq\chi_{-j}\boxtimes\Sym^{2j}V$, and the
second is $\mathcal S_j^-\simeq\chi_j\boxtimes\Sym^{2j}V^*$.
Multiplication by $a^q$ gives degree $j+q$ and circle weight $-j-2q$;
multiplication by $b^q$ on the dual side gives weight $j+2q$.
The fibre product in \eqref{eq:ring6} retains one common constant.
This proves the graded $K_1\times K_2$-decomposition and completes
Theorem~\ref{thm:quotient} for $E_6$.

\section{The linear section in type \texorpdfstring{$E_7$}{E7}}\label{sec:e7}
Here we compute $Q=R/(\fp_2)$ for the long-root
$\SL_2\times\Spin_{12}$ pair in $E_7$. The geometry of four-planes
in a direct sum gives two Segre cones; their equations and graded
coordinate rings prove Theorem~\ref{thm:quotient} in this case.

\subsection{The two Segre cones}
Let $U=V\oplus W$, with $\dim V=\dim W=4$. A point
$L\in\Gr(4,U)$ is a four-dimensional subspace of $U$.
The vector decomposition
\[
 \Lambda^4U=\Lambda^4V\oplus(\Lambda^3V\otimes W)
 \oplus(\Lambda^2V\otimes\Lambda^2W)
 \oplus(V\otimes\Lambda^3W)\oplus\Lambda^4W
\]
has weights $2,1,0,-1,-2$; its coordinate functions have the opposite
weights. Here $\fp_2=\Lambda^2V\otimes\Lambda^2W$ is precisely
the middle term. Under the Killing-form identification with linear
functions, the equations $\fp_2=0$ set all middle Pl\"ucker coordinates
to zero. Thus the middle component of a Pl\"ucker vector for $L$ must
vanish. We first describe the resulting points geometrically; afterwards
we check the scheme structure with equations.

For a four-plane $L\subset V\oplus W$, put
$r=\dim(L\cap V)$ and $s=\dim(L\cap W)$. Choose bases of these
two intersections and complete them to a basis of $L$ by vectors
$v_i+w_i$. The classes of the $v_i$ in $V/(L\cap V)$ are linearly
independent, as are the classes of the $w_i$ in $W/(L\cap W)$.
Expanding the wedge of this basis, its nonzero components have exactly
the $V$-degrees $r,r+1,\ldots,4-s$: the terms of each such degree
are distinct basis wedges and cannot cancel. In particular the middle
component is nonzero exactly when $r\leq2$ and $s\leq2$.
Its vanishing therefore gives
\[
 X_+=\{L:\dim(L\cap V)\geq3\},\qquad
 X_-=\{L:\dim(L\cap W)\geq3\}
 \quad\text{in }\Gr(4,U).
\]
These loci are disjoint, since a four-plane cannot contain both a
three-plane in $V$ and a three-plane in $W$.

Now let $L\in X_+$, with $L\ne V$. Then $A=L\cap V$ has
dimension three, and $L/A$ is a line in
\[
 U/A=(V/A)\oplus W.
\]
Its projection to $W$ is nonzero, so its image is a line $D\subset W$.
We have therefore associated to $L$ a pair
\[
 (A,D)\in\Gr(3,V)\times\Gr(1,W)
       =\PP(V^*)\times\PP(W).
\]
Conversely, fix such a pair. The four-planes containing $A$ and contained
in $V\oplus D$ are parameterized by
$\PP((V/A)\oplus D)\simeq\PP^1$. One of them is $V$; all the
others have intersection $A$ with $V$ and image $D$ in $W$.
Thus fixing $(A,D)$ gives a whole line of four-planes, with $V$ as
one endpoint and $A\oplus D$ as the other distinguished point.

This is also a projective line in the Pl\"ucker embedding. Indeed,
choose a basis $u_1,u_2,u_3$ of $A$, a vector $v\in V\setminus A$,
and a nonzero $w\in D$. The four-plane
$A+\C(\lambda v+\mu w)$ has Pl\"ucker point
\[
 [\lambda\,u_1\wedge u_2\wedge u_3\wedge v
   +\mu\,u_1\wedge u_2\wedge u_3\wedge w],
 \qquad [\lambda:\mu]\in\PP^1.
\]
The first endpoint is the fixed point $[\Lambda^4V]$, representing $V$.
The second endpoint, representing $A\oplus D$, varies over the Segre
embedding of $\Gr(3,V)\times\Gr(1,W)$ in
$\PP(\Lambda^3V\otimes W)$. Hence $X_+$ is the projective cone
over this product, with vertex $V$. Projection from the vertex is
exactly $L\mapsto(A,D)$. Interchanging $V$ and $W$ gives $X_-$,
the cone over $\PP(V)\times\PP(W^*)$ with vertex $W$.

We now verify that these two cones give the entire linear section
scheme-theoretically, with no additional nilpotents.

Index bases of $V,W$ by $1,\ldots,4$ and $5,\ldots,8$.
Let $I$ be the Pl\"ucker ideal in $\C[p_J]$, and let $I_2$ be its
homogeneous degree-two part, the vector space of Pl\"ucker quadrics.
These quadrics generate $I$. Set to zero the coordinates with two
indices in each block, and denote the image of $I_2$ by
$\overline I_2$. The surviving coordinates are $a=p_{1234}$,
$b=p_{5678}$, and two $4\times4$ arrays: $x$ has three $V$ indices
and one $W$ index, while $y$ has one $V$ index and three $W$ indices.
Thus $\overline I_2$ is a space of quadratic polynomials in
$S=\C[a,b,x,y]$, and $Q=S/(\overline I_2)$, where parentheses
mean the ideal generated by this space.

Write $\mathcal X=\operatorname{span}\{x_{ij}\}$ and
$\mathcal Y=\operatorname{span}\{y_{ij}\}$. After the cofactor
sign conventions, these are the $\SU(V)\times\SU(W)$-modules
$V\otimes W^*$ and $V^*\otimes W$, respectively; $a,b$ are invariant.
Let $\overline I_2(x)$ be the vector space spanned by the quadratic minors
\[
 x_{ij}x_{kl}-x_{il}x_{kj}\qquad(i<k,\ j<l),
\]
and define $\overline I_2(y)$ similarly. In particular,
\[
 \overline I_2(x)\simeq\Lambda^2V\otimes\Lambda^2W^*,\qquad
 \overline I_2(y)\simeq\Lambda^2V^*\otimes\Lambda^2W,
\]
each of dimension $6\cdot6=36$.
The decomposition of the specialized quadratic relation space is
\begin{equation}\label{eq:ideal7}
 \overline I_2=\overline I_2(x)\oplus a\mathcal Y
 \oplus\bigl(\C ab\oplus\mathcal X\mathcal Y\bigr)
 \oplus b\mathcal X\oplus\overline I_2(y).
\end{equation}
Here $\mathcal X\mathcal Y$ denotes the span of all products
$x_{ij}y_{kl}$, naturally isomorphic to $\mathcal X\otimes\mathcal Y$.
All summands in \eqref{eq:ideal7} are spaces of quadratic generators.
They are grouped by total $V$-degree, rather than into irreducible
representations: the middle summand is reducible.

To verify the decomposition, substitute the middle-coordinate equations
in the generating Pl\"ucker quadrics
\begin{equation}\label{eq:plucker}
 \sum_{r=1}^{5}(-1)^{r-1}
 p_{i_1i_2i_3j_r}\,p_{j_1\cdots\widehat{j_r}\cdots j_5}=0.
\end{equation}
Repeated indices give zero and permutation of indices changes signs.
Grouping by the total number of $V$ indices gives the following summands:
\begin{center}
\begin{tabular}{c c c c}
\toprule
total $V$-degree & quadratic generators & $\SU(V)\times\SU(W)$-module & dimension \\
\midrule
$6$ & $\overline I_2(x)$ & $\Lambda^2V\otimes\Lambda^2W^*$ & $36$ \\
$5$ & $a\mathcal Y$ & $V^*\otimes W$ & $16$ \\
$4$ & $\C ab\oplus\mathcal X\mathcal Y$
 & $1\oplus\bigl(\operatorname{End}V\otimes\operatorname{End}W\bigr)$ & $257$ \\
$3$ & $b\mathcal X$ & $V\otimes W^*$ & $16$ \\
$2$ & $\overline I_2(y)$ & $\Lambda^2V^*\otimes\Lambda^2W$ & $36$ \\
\bottomrule
\end{tabular}
\end{center}
The dimensions follow immediately from the displayed modules:
$6\cdot6$, $4\cdot4$, $1+16\cdot16$, $4\cdot4$, and $6\cdot6$.
We verify directly that these spaces are exactly the specialized
relations.

The minors of $x,y$ and the products $ay_{ij},bx_{ij}$ follow at once
from \eqref{eq:plucker} after setting the middle coordinates to zero.
For the mixed products, write $v_1,\ldots,v_4$ and $w_1,\ldots,w_4$
for the two indexed bases, with index sets
$\mathcal V=\{v_1,\ldots,v_4\}$ and
$\mathcal W=\{w_1,\ldots,w_4\}$. Choose the cofactor signs
\[
 x_{ij}=(-1)^{i-1}p_{v_1\cdots\widehat{v_i}\cdots v_4w_j},
 \qquad
 y_{ij}=(-1)^{j-1}p_{v_iw_1\cdots\widehat{w_j}\cdots w_4}.
\]
If $i\ne k$, choose a third $V$ index $m$. In \eqref{eq:plucker}
take the three-index list to be
$(\mathcal V\setminus\{v_i,v_m\})\cup\{w_j\}$ and the five-index list to be
$\{v_m,v_k\}\cup(\mathcal W\setminus\{w_l\})$.
Only the term $x_{ij}y_{kl}$ survives, so this product vanishes.
Interchanging the blocks gives the same conclusion when $j\ne l$.
For the remaining products, fix $j$ and put $c_i=x_{ij}y_{ij}$.
Using the same three-index list and the five-index list
$\{v_m,v_i\}\cup(\mathcal W\setminus\{w_j\})$ gives
$c_i+c_m=0$ for every $i\ne m$.
Taking three distinct indices yields $c_i=0$. Thus all $xy$ products
are relations. Finally, the three-index list $\mathcal V\setminus\{v_i\}$
and the five-index list $\{v_i\}\cup\mathcal W$ give $ab$ plus a linear
combination of these products, so $ab$ is a relation as well.

Conversely, the proposed cones consist of decomposable four-vectors,
so every specialized Pl\"ucker quadric vanishes on them.
Their quadratic vanishing space is precisely the right side of
\eqref{eq:ideal7}: within each coordinate block the relations are the
minors, and between the blocks they are all mixed products.
This proves \eqref{eq:ideal7}. Since $I$ is generated by $I_2$,
its specialization is generated by $\overline I_2$.

In particular, the relations of total $V$-degrees $5,4,3$ are exactly
\[
 ay_{ij}=0,\qquad ab=0,\quad x_{ij}y_{kl}=0,\qquad bx_{ij}=0.
\]
Together they generate the product ideal $(a,x)(b,y)$.
Thus a point has either $(a,x)=0$ or $(b,y)=0$; the remaining minors
cut out the corresponding rank-one cone. The two affine components
meet only at the origin, and their projectivizations are disjoint.

Let $(\overline I_2(x))\subset\C[x]$ and
$(\overline I_2(y))\subset\C[y]$ be the ideals generated by the two
spaces of minors. Adjoining the vertex coordinates and imposing the
mixed quadratic generators gives
\begin{equation}\label{eq:ring7}
 Q=\bigl(\C[x]/(\overline I_2(x))\bigr)[a]\times_\C
       \bigl(\C[y]/(\overline I_2(y))\bigr)[b].
\end{equation}
Both factors are polynomial extensions of rank-one matrix rings, hence
reduced. This proves the scheme assertion and the cone description.

The substitution $x_{ij}=s_it_j$ identifies the degree-$j$ piece
of the first matrix ring with the polynomials of bidegree $(j,j)$:
\[
 \mathcal S_j^+=\bigl(\C[x]/(\overline I_2(x))\bigr)_j
 \simeq\chi_{-j}\boxtimes
       \bigl(\Sym^jV\otimes\Sym^jW^*\bigr).
\]
The second ring gives the dual $K_2$-type with circle weight $j$.
Adjoining $a,b$, of degree one and weights $-2,2$, gives the terms
$a^q\mathcal S_j^+$ and $b^q\mathcal S_j^-$ in degree $j+q$,
with weights $-j-2q$ and $j+2q$. The fibre product in
\eqref{eq:ring7} retains one common constant, proving
Theorem~\ref{thm:quotient} for $E_7$.

\section{The linear section in type \texorpdfstring{$E_8$}{E8}}\label{sec:e8}
This section computes $Q=R/(\fp_2)$ for the long-root
$\SL_2\times E_7$ pair in $E_8$. We identify two cones over
Grassmannians inside the spinor variety, verify their scheme structure,
and read off the graded coordinate rings to prove
Theorem~\ref{thm:quotient}.

\subsection{Pure spinors and the two Grassmannian cones}
Write the orthogonal space as $U=V\oplus V^*$, $\dim V=8$,
with the hyperbolic symmetric form
\[
 \langle(v,\phi),(w,\psi)\rangle=\phi(w)+\psi(v).
\]
The variety $\OGr^+(8,U)$ parameterizes one of the two components
of maximal isotropic eight-planes; we choose the component containing
$V$. It also contains $V^*$, since $\dim V=8$ is even.
We write $L\in\OGr^+(8,U)$ for such an eight-dimensional maximal
isotropic subspace of $U$. We use the spinor embedding of $\OGr^+(8,U)$.
The subgroup of $\SO(U)$ acting by
$(v,\phi)\mapsto(gv,(g^{-1})^*\phi)$, with $g\in\operatorname{GL}(V)$,
is a Levi subgroup of the parabolic stabilizing $V$.
By the \emph{covering Levi} we mean its inverse image in $\Spin(U)$.
Explicitly, this group is
\[
 \widetilde{\operatorname{GL}}(V)
   =\{(g,z)\in\operatorname{GL}(V)\times\C^\times:z^2=\det g\}.
\]
Projection $(g,z)\mapsto g$ is a double cover with kernel
$\{(1,1),(1,-1)\}$. Its compact form is the double cover of
$\operatorname{U}(V)$ obtained by requiring $g$ unitary and $|z|=1$.
Here $V$ has dimension eight, whereas
$U=V\oplus V^*$ has dimension sixteen.
The one-dimensional representation $\delta(g,z)=z$ is an actual
character of this cover and satisfies $\delta^2=\det V$.
The exterior-algebra construction gives the ambient Cartan space
$\fp=S^+=\delta\otimes\Lambda^{\rm even}V^*$.
On the lifted copy $g\mapsto(g,1)$ of $\SU(V)$, the character
$\delta$ is trivial; a volume form
identifies $\Lambda^6V^*\simeq\Lambda^2V$ and
$\Lambda^4V^*\simeq\Lambda^4V$. Thus, as an $\SU(V)$-module,
\[
 \fp\simeq\C\oplus\Lambda^2V^*
       \oplus\underbrace{\Lambda^4V}_{\fp_2}
       \oplus\Lambda^2V\oplus\C.
\]
These five summands come from exterior degrees $0,2,4,6,8$ in the
spinor model. Their dimensions are $1,28,70,28,1$, and their circle
weights are $2,1,0,-1,-2$: the factor $\delta$ has weight $2$,
while $V^*$ has weight $-1/2$ on the cover.
The half-spin action descends to the semispin quotient used for $K$
in Section~\ref{sec:setup}.

Dually, the degree-one coordinate space is
$R_1=(S^+)^*=\delta^{-1}\otimes\Lambda^{\rm even}V$,
whose five pieces have circle weights $-2,-1,0,1,2$.
The middle piece is the $70$-dimensional Cartan space $\fp_2$ of the
centralizer, identified with linear functions by the Killing form.
Thus $\fp_2=0$ sets all exterior-degree-four spinor coordinates to zero.
Let $a,x,y,b$ denote the remaining coordinate pieces, of exterior
degrees $0,2,6,8$, respectively. Using the fixed volume form, identify
the $y$ coordinates with complementary two-index coordinates on the
dual side.

We first examine the open chart $a\ne0$, consisting of eight-planes
transverse to $V^*$. Such a plane is the graph of a map
$B:V\to V^*$. Isotropy says
$B(v)(w)+B(w)(v)=0$, so $B$ is an alternating form, an element of
$\Lambda^2V^*$. The spinor coordinates on this chart are its
principal Pfaffians; equivalently its spinor line is represented by
\[
 \delta\exp(B)=\delta\left(1+B+\frac{B\wedge B}{2!}
                  +\frac{B^{\wedge3}}{3!}+\frac{B^{\wedge4}}{4!}\right).
\]
We use the usual Pfaffian chart conventions, as in
\cite[Section~2]{Rincon}. Killing the middle coordinates says
$B\wedge B=0$. An alternating form has this property exactly when
its rank is at most two. Thus $B=0$, giving $L=V$, or $B$ is
a decomposable two-form with six-dimensional kernel. Since the graph
meets $V$ in $\ker B$, this calculation leads to the locus $X_+$
below; the opposite chart gives $X_-$:
\[
 X_+=\{L\in\OGr^+(8,U):\dim(L\cap V)\geq6\},\qquad
 X_-=\{L\in\OGr^+(8,U):\dim(L\cap V^*)\geq6\}.
\]
They are disjoint, since $6+6>8$. We next describe each entire locus,
including the boundary of its graph chart. The equations in the next
subsection will show that these are exactly the two components of
the linear section, with no further components or nilpotents.

Take $L\in X_+$ with $L\ne V$. In our chosen orthogonal component
$\dim(L\cap V)$ is even, so it must be six. Set
\[
 A=L\cap V,\qquad Q_A=V/A,\qquad
 D=\operatorname{Ann}(A)=Q_A^*\subset V^*.
\]
Both $Q_A$ and $D$ have dimension two. Isotropy forces $L$ to lie
in $A^\perp=V\oplus D$, so $L/A$ is a maximal isotropic two-plane
in the four-dimensional orthogonal space
\[
 A^\perp/A=Q_A\oplus Q_A^*.
\]
It belongs to the same component as $Q_A$.

For fixed $A$, this component is a projective line. To see this
directly, its two-planes transverse to $Q_A^*$ are graphs of alternating
maps $Q_A\to Q_A^*$. These maps form the one-dimensional space
$\Lambda^2Q_A^*=\det Q_A^*$. The only remaining two-plane in
this component is $Q_A^*$ itself: intersection with $Q_A^*$ has
even dimension, so a nontransverse plane must equal it.
Thus the affine line of graphs acquires one point at infinity, giving
\[
 \OGr^+(2,Q_A\oplus Q_A^*)
        =\PP(\C\oplus\det Q_A^*)\simeq\PP^1.
\]
The zero map gives $Q_A$, hence $L=V$; the point at infinity gives
$Q_A^*$, hence $L=A\oplus D$. Every other member has intersection
exactly $A$ with $V$.

This projective line is a line in the spinor embedding as well.
Choose $0\ne\eta\in\det D\subset\Lambda^2V^*$.
Since $\eta\wedge\eta=0$, the graphs of multiples of $\eta$
have spinor points $[\delta(1+c\eta)]$. Their closure is the line
\[
 [\delta(\lambda+\mu\eta)],\qquad[\lambda:\mu]\in\PP^1,
\]
joining the vertex $[\delta]$, representing $V$, to
$[\delta\eta]$, representing $A\oplus D$. As $A$ varies, its
annihilator $D$ varies through $\Gr(2,V^*)$, and $[\eta]=[\det D]$
is precisely its Pl\"ucker point. Consequently $X_+$ is the
projective cone over the Pl\"ucker-embedded $\Gr(2,V^*)$, with
vertex $V$. Projection from the vertex sends $L\ne V$ to
$D=\operatorname{Ann}(L\cap V)$. The description of $X_-$ is
dual, with base $\Gr(2,V)$ and vertex $V^*$.

\subsection{The specialized quadratic generators}
Put $N=\{1,\ldots,8\}$ and label spinor coordinates $p_I$ by even
subsets $I\subset N$. Let $I$ be the pure-spinor ideal and $I_2$
its homogeneous degree-two part. This space is spanned by the Wick
relations and generates $I$ \cite[Section~2]{Rincon}.
For odd subsets $S,T$, write
$S\mathbin\triangle T=\{i_1<\cdots<i_m\}$; the relation is
\begin{equation}\label{eq:wick}
 \sum_{r=1}^m(-1)^{r-1}
 p_{S\mathbin\triangle\{i_r\}}
 p_{T\mathbin\triangle\{i_r\}}=0.
\end{equation}
Set all $p_I$ with $|I|=4$ to zero, and denote the image of $I_2$
by $\overline I_2$. The surviving coordinates are $a=p_\varnothing$,
$b=p_N$, and two skew arrays, which we normalize by
\[
 x_{ij}=p_{\{i,j\}},\qquad
 y_{ij}=(-1)^{i+j-3}p_{N\setminus\{i,j\}}\quad(i<j),
 \qquad x_{ji}=-x_{ij},\quad y_{ji}=-y_{ij}.
\]
The sign on $y$ is the complement sign for the fixed volume form.
As $\SU(V)$-modules, the coordinate spans are
$\mathcal X\simeq\Lambda^2V$ and $\mathcal Y\simeq\Lambda^2V^*$,
each of dimension $28$. The full specialized ideal in
$\C[a,b,x,y]$ is $(\overline I_2)$.

Write $\overline I_2(x)$ for the space spanned by the quadratic
Pfaffians
\[
 x_{ij}x_{kl}-x_{ik}x_{jl}+x_{il}x_{jk},\qquad i<j<k<l,
\]
and define $\overline I_2(y)$ similarly. These are the Pl\"ucker
quadrics for $\Gr(2,V^*)$ and $\Gr(2,V)$, respectively.
Their module structures are
\[
 \overline I_2(x)\simeq\Lambda^4V,\qquad
 \overline I_2(y)\simeq\Lambda^4V^*\simeq\Lambda^4V.
\]
Indeed, the Pfaffian map identifies the first space with the
$\Lambda^4V$ summand in
$\Sym^2(\Lambda^2V)=\bS_{(2,2)}V\oplus\Lambda^4V$;
the second is dual. Each has dimension $\binom84=70$.

The full space of specialized quadratic generators is
\begin{equation}\label{eq:ideal8}
 \overline I_2=\overline I_2(x)\oplus a\mathcal Y
 \oplus\bigl(\C ab\oplus\mathcal X\mathcal Y\bigr)
 \oplus b\mathcal X\oplus\overline I_2(y).
\end{equation}
Here $\mathcal X\mathcal Y\simeq\Lambda^2V\otimes\Lambda^2V^*$
is the span of all $xy$ products. With $V_\lambda$ denoting the
irreducible $\SU(8)$-module of highest weight $\lambda$, it decomposes as
\[
 \mathcal X\mathcal Y
 =\C\oplus V_{\omega_1+\omega_7}\oplus V_{\omega_2+\omega_6},
 \qquad \dim=1+63+720.
\]
Consequently the summands of \eqref{eq:ideal8} are as follows:
\begin{center}
\begin{tabular}{c c c c}
\toprule
circle weight & generator space & $\SU(8)$-module & dimension \\
\midrule
$-2$ & $\overline I_2(x)$ & $V_{\omega_4}$ & $70$ \\
$-1$ & $a\mathcal Y$ & $V_{\omega_6}$ & $28$ \\
$0$ & $\C ab\oplus\mathcal X\mathcal Y$
 & $\C^{\oplus2}\oplus V_{\omega_1+\omega_7}\oplus V_{\omega_2+\omega_6}$ & $785$ \\
$1$ & $b\mathcal X$ & $V_{\omega_2}$ & $28$ \\
$2$ & $\overline I_2(y)$ & $V_{\omega_4}$ & $70$ \\
\bottomrule
\end{tabular}
\end{center}
The two scalar summands at weight zero are $\C ab$ and the invariant
contraction in $\mathcal X\mathcal Y$. The dimensions follow from the
displayed modules and sum to $981$.

We prove \eqref{eq:ideal8} directly from \eqref{eq:wick}.
Taking $S=\{i\}$ and $T=\{j,k,l\}$, with four distinct indices,
gives a Pfaffian in $x$, since the term involving $p_{\{i,j,k,l\}}$
vanishes. For a six-element set $C$, taking $S=\{i\}$ and
$T=C\setminus\{i\}$ with $i\in C$ leaves only $ap_C$;
all other terms contain a four-index coordinate. This gives $a\mathcal Y$.
Complementing all index sets gives the dual Pfaffians and $b\mathcal X$.

For the mixed products there are three cases. If $i,j,k,l$ are distinct,
choose $r$ outside these four indices and take
$S=\{i,j,r\}$ and $T=N\setminus\{k,l,r\}$.
Only the term proportional to $x_{ij}y_{kl}$ survives, so it is a relation.
If the pairs share one index, fix distinct $i,l$ and put
$q_j=x_{ij}y_{jl}$ for $j\notin\{i,l\}$.
The choices $S=\{i,j,k\}$ and $T=N\setminus\{j,k,l\}$ give
$q_j+q_k=0$ for distinct $j,k$ outside $\{i,l\}$.
Using three such indices yields $q_j=0$.
Finally, put $c_{ij}=x_{ij}y_{ij}=c_{ji}$.
Taking $S=\{i,j,k\}$ and $T=N\setminus S$ gives
\[
 c_{ij}+c_{ik}+c_{jk}=0.
\]
Comparing the four triangle equations on four distinct indices yields
$c_{ij}=c_{kl}$ for disjoint pairs. Any two pairs have a third pair
disjoint from both, so all $c_{ij}$ are equal. The triangle equation
then forces them all to vanish. Thus every $xy$ product is a relation.
Taking $S=\{i\}$ and $T=N\setminus\{i\}$ now gives $ab$ plus
mixed products already shown to vanish, proving that $ab$ is a relation.

Conversely, the cone in the $(a,x)$ coordinates consists of pure
spinors: for $a\ne0$ it is represented by $a\exp(B)$ with $B$
decomposable, and the remaining points lie in its closure.
The other cone is dual. Hence every specialized Wick quadric vanishes
on these two cones. Their quadratic vanishing space is the right side
of \eqref{eq:ideal8}: the Pfaffians within each block and all products
between the blocks. This proves equality of the generator spaces.
Since $I=(I_2)$, these quadrics generate the entire specialized ideal.

The generators at circle weights $-1,0,1$ generate $(a,x)(b,y)$.
They separate the two components; the remaining generators impose
decomposability on each skew array. Writing parentheses for the ideals
generated by the pure quadratic spaces, we obtain
\begin{equation}\label{eq:ring8}
 Q=\bigl(\C[x]/(\overline I_2(x))\bigr)[a]\times_\C
       \bigl(\C[y]/(\overline I_2(y))\bigr)[b].
\end{equation}
Here $\C[x]=\C[x_{ij}:i<j]$, and similarly for $y$.
The two matrix rings are the reduced homogeneous coordinate rings
of the Pl\"ucker-embedded $\Gr(2,V^*)$ and $\Gr(2,V)$.
Adjoining $a,b$ gives their projective cones with vertices $V,V^*$.
Both maps in the fibre product evaluate at the origins. Thus the two
affine components meet only at the origin, while their projectivizations
are disjoint. This proves the scheme assertion.

The standard graded decomposition of the Pl\"ucker coordinate ring
of $\Gr(2,V^*)$ gives
\[
 \mathcal S_j^+=\bigl(\C[x]/(\overline I_2(x))\bigr)_j
       \simeq\chi_{-j}\boxtimes\bS_{(j,j)}V.
\]
Here each $x$-coordinate has degree one and circle weight $-1$.
The second ring gives $(\bS_{(j,j)}V)^*$ with weight $j$.
Thus $a^q\mathcal S_j^+$ and $b^q\mathcal S_j^-$ have degree
$j+q$ and weights $-j-2q$ and $j+2q$, respectively.
The fibre product in \eqref{eq:ring8} retains one common constant.
This proves Theorem~\ref{thm:quotient} for $E_8$.

\bigskip
{\small
\noindent Mathematics Department, University of Utah, Salt Lake City, UT 84112\par
\noindent\textit{Email address:}
\href{mailto:gordan.savin@utah.edu}{\texttt{gordan.savin@utah.edu}}\par
}
\end{document}